\documentclass{amsart}
\usepackage{amssymb}
\usepackage{graphicx}
\usepackage[all]{xy}
\newtheorem*{theorem*}{Theorem}
\newtheorem{theorem}{Theorem}[section]
\newtheorem{lemma}[theorem]{Lemma}

\newtheorem{corollary}[theorem]{Corollary}

\theoremstyle{definition}
\newtheorem{definition}[theorem]{Definition}

\theoremstyle{remark}
\newtheorem{remark}[theorem]{Remark}

\numberwithin{equation}{section}

\newcommand{\firef}[1]{Figure~{\rm\ref{#1}}}

\newcommand{\leref}[1]{Lemma~{\rm\ref{#1}}}

\newcommand{\deref}[1]{Definition~{\rm\ref{#1}}}

\newcommand{\fig}[1]
{\raisebox{-0.5\height}%
{\includegraphics{figures/#1}}}
\newcommand{\st}{\; | \;}                     
\newcommand{\ttt}{\otimes}                    

\newcommand{\del}{\partial}
\newcommand{\<}{\langle}
\renewcommand{\>}{\rangle}
\newcommand{\isoto}{\xrightarrow{\sim}}       
\newcommand{\xxto}{\xrightarrow}              

\newcommand{\tbox}{\boxtimes}
\newcommand{\dotimes}{\otimes\dots\otimes} 
\renewcommand{\1}{\mathbf{1}}      
\renewcommand{\Vec}{\mathcal{V}ec} 
\newcommand{\CC}{{\mathcal C}}     
\newcommand{\G}{{\mathcal G}}      
\newcommand{\C}{\mathbb{C}}       
\newcommand{\Z}{\mathbb{Z}}       
\newcommand{\bg}{\mathbf{g}}      
\newcommand{\bh}{\mathbf{h}}      
\newcommand{\bx}{\mathbf{x}}      

\newcommand{\tpa}{\tilde{p}_a}

\newcommand{\tsi}{\widetilde \Sigma}
\newcommand{\RR}{\mathcal{R}}       

\newcommand{\M}{M(\widetilde \Sigma,\Sigma)} 
\newcommand{\al}{\alpha}
\newcommand{\be}{\beta}

\newcommand{\ph}{\varphi}

\newcommand{\si}{\sigma}
\newcommand{\Si}{\Sigma}

\renewcommand{\th}{\theta}

\DeclareMathOperator{\id}{id}

\DeclareMathOperator{\Hom}{Hom}

\DeclareMathOperator{\End}{End}

\begin{document}

\title{On $G$--modular functor}

\author{Alexander Kirillov, Jr. \and Tanvir Prince}
   \address{Department of Mathematics, SUNY at Stony Brook,
            Stony Brook, NY 11794, USA}
    \email{kirillov@math.sunysb.edu}
    \urladdr{http://www.math.sunysb.edu/\textasciitilde kirillov/}
    \address{Graduate student. Advisor:Alexander Kirillov, Jr.
Department of Mathematics, SUNY at Stony Brook,
            Stony Brook, NY 11794, USA}
    \email{prince@math.sunysb.edu}

\date{\today}

\maketitle
\begin{abstract}
It is a known result that given a $\CC$ extended modular functor where
$\CC$ is a semisimple abelian category, we can find a structure of
weakly rigid fusion category on $\CC$. Also if we have a structure of a
weakly rigid fusion category on $\CC$ then from this we can define a
$\CC$ extended modular functor. In this paper, we extend this notion of
modular functor and fusion category to what we called $G$ equivariant
modular functor and $G$ equivariant fusion category where $G$ is a
finite group. Then we establish a similar correspondence between $G$
equivariant modular functor and $G$ equivariant fusion category.
\end{abstract}

\tableofcontents

\section{Introduction}
This paper is a continuation of \cite{equimf1} and \cite{TA}. Its goal
is to develop
a formalism of $G$-equivariant modular functors, which would provide a
suitable algebraic formalism for orbifold models in conformal field
theory, much as usual modular functors can be used for describing
various structures appearing in usual conformal field theory. We will
also establish a relation of this approach to the theory of
$G$-equivariant fusion categories as defined in \cite{equimf1}
(following earlier work of Turaev \cite{T2}). The main result of this
paper is that the notion of $G$-equivariant modular functor (in genus
zero) and a structure of a $G$-equivariant fusion category are
essentially equivalent; precise statement is given in the main theorem
of this paper. 

Our approach only discusses topological setting, in which the main
objects are oriented surfaces with boundary (or, in $G$-equivariant
case, $G$-covers of such curves).  Complex-analytic analog,
which uses the language of flat connections on the moduli spaces of
curves,  will be discussed in a subsequent paper.

It should be noted that some of the results here are parallel to the
results of \cite{T2}. However, unlike Turaev, our approach is not
based on 3D TQFT, in which the main technical tool is presentation of
3-manifolds via surgery and using Kirby moves. Instead, we follow the
approach suggested (for  non-$G$-equivariant case) by Moore and
Seiberg,  presenting a surface as a result of gluing of ``standard''
spheres with holes, and then writing  the generators and relations in
the groupoid of all such presentations.

\section{$G$-equivariant fusion category}

Throughout this paper, $G$ is  a fixed finite group.

\begin{figure}[[ht]
\fig{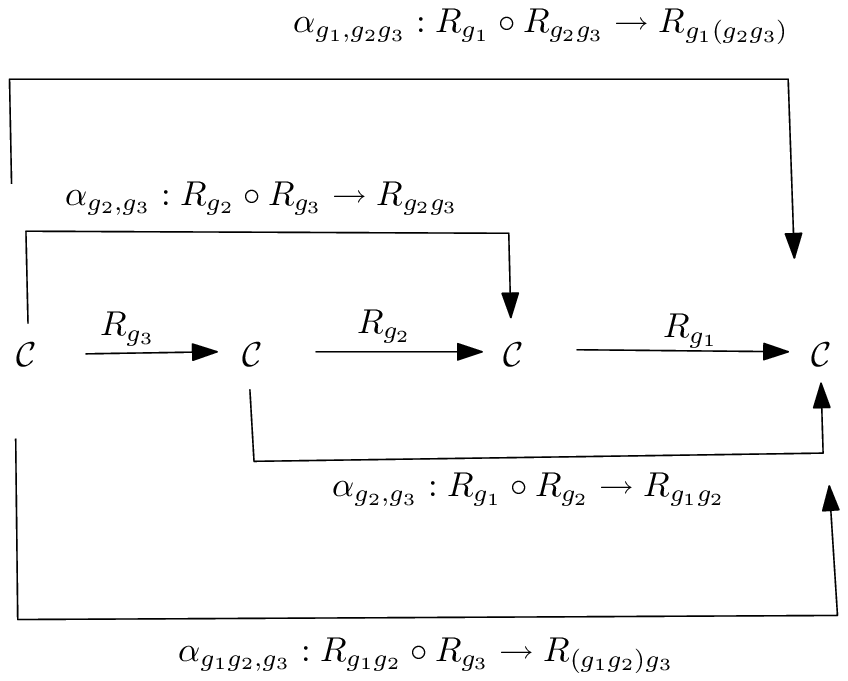}
\caption{Associativity of $G$ action}\label{assg}
\end{figure}

\begin{definition}\label{geqcat}
  A {\em $G$-equvariant category} is an abelian category
  $\CC$ with the following additional structure:
  \begin{description}
    \item[$G$-\textbf{grading}] Decomposition
      $$\CC = \bigoplus_{g \in G}\CC_{g}$$
      where each $\CC_{g}$ is a full subcategory of $\CC$. We will
      call objects $V \in \CC_{g}$ ``$g$ twisted''. In particular,
      objects $V \in \CC_{1}$ will be called ``neutral''. In physical
      literature, the subcategory $\CC_{1}$ is usually called the
    ``untwisted sector''.
      
    \item[\textbf{Action of $G$}] For each $g \in G$, we are given a
      functor $R_{g}\colon \CC \to \CC$, and for each pair $g,h \in
      G$, a  functorial isomorphism $\al_{g,h}\colon R_{g} \circ
      R_{h} \to R_{gh}$. These functorial isomorphisms must satisfy
      the following conditions:
      \begin{enumerate}
        \item $R_{1} = \id$
        \item $R_{g}(\CC_{h})\subseteq \CC_{ghg^{-1}}$. This means the
          action of $G$ respects the grading.
        \item $\al_{{g_1g_2},{g_3}} \circ \al_{g_1,g_2} =
          \al_{{g_1},{g_2g_3}} \circ \al_{g_2,g_3}$. Here both
          sides are functorial isomorphism from $R_{g_1}R_{g_2}R_{g_3}
          \to R_{g_1g_2g_3}$ (see  \firef{assg}). This might be thought
            of as the  associativity of $G$ action. 
      \end{enumerate}
     We will frequenlty use notation ${}^gV$ for $R_g(V)$. 
  \end{description}
\end{definition}
This definition was introduced by Turaev \cite{T2} under a differnt
name. 

\begin{definition}\label{geqfucat}
  A {\em $G$-equivariant fusion category} is a semisimple $G$
  equivariant  abelian category over the base field $\C$ with the
  following  additional structure:
  \begin{enumerate}
    \item Structure of a  monoidal category such that\\
      $\1$ is a simple object\\
      for any simple object $V_i$, $\End_{\CC}(V_i)=\C$\\
      $R_g$ is a tensor functor\\
      For $X \in \CC_{g}$ and $Y \in \CC_{h}$,  $X \otimes Y \in
      \CC_{gh}$. 
    \item Structure of rigidity:  each object $V$ has a right
      dual $V^{*}$, with evaluation and coevaluation map, 
      $e_V\colon V^* \otimes V \to \1$ and 
      $i_V\colon \1 \to V \otimes V^*$ which satisfy the following
      rigidity conditions:\\
      $(\id \otimes e_V)\circ (i_V \otimes \id)$ is the identity map on
        $V$\\
      $(e_V \otimes \id) \circ (\id \otimes i_V)$ is the identity map on
          $V^*$\\
      For detailed discussion on this, see e.g. \cite{BK}.
    
    \item Functorial isomorphism $\delta_{V}:V \to V^{**}$, satisfying
      the same compatibility conditions as in the absence of $G$. For
      detailed discussion see \cite{BK}. Here, for readers convenience,
      we list these conditions:\\
        $\delta_{V \otimes W} = \delta_{V} \otimes \delta_{W}$\\
        $\delta_{\1} = \id$\\
        $\delta_{V^*} = (\delta_V^*)^{-1}$\\
        In addition to the above reations, we also require\\
        $R_g(\delta_V) = \delta_{R_g(V)}$

    \item A collection of functorial isomorphisms, $R_{V,W}\colon  V
        \otimes W \to {}^gW \otimes V$ for every $V \in \CC_{g}$
        and $W \in \CC_{h}$. This is similar to the braiding in
        the $G = \{1\}$ case  but  with the addtion of the twist. These
        functorial  isomorphisms must satisfy an analog of two
        hexagon axioms.
        The first hexagon axiom is shown in \firef{hexagon_fusion};
        the other one is similar but with $R$ replaced by $R^{-1}$. 

            \begin{figure}[ht]
            $$\xymatrix{
                        &(^hV\otimes U)\otimes W\ar[r]^{\al} & 
                                      {}^hV\otimes (U\otimes W)
                                  \ar[dr]^{\id\otimes R_{U,W}}\\
            (U\otimes V)\otimes W
               \ar[ur]^{R_{U,V}\otimes\id}\ar[dr]^\al&&& 
                                      {}^hV\otimes ({}^hW\otimes U)\\
                      &U\otimes (V\otimes W) \ar[r]^{R_{U,V\otimes W}}&
                                      {}^h(V\otimes W)\otimes U
                                      \ar[ur]_\al 
            }$$
            \caption{Hexagon axiom ($U\in \CC_h$)}\label{hexagon_fusion}
            \end{figure}
  \end{enumerate}
\end{definition}

The definition immediately implies that $\1 \in \CC_{1}$ and if $V \in
\CC_{g}$ then $V^{*} \in \CC_{g^{-1}}$.
Also, since in a rigid monoidal category the unit object and dual is
unique up to a unique isomorphism, we have canonical identification
\begin{align*}
  {}^g\1& = \1 \\
  ({}^gV)^*& = {}^g(V^*).
\end{align*}

\begin{remark}
  From now on we will refer to the isomorphisms described above:
  associativity, unit, braiding, rigidity,  $\delta$ morphism,
  $G$-action, and their  compositions as canonical morphisms. We will
  omit these canonical  morphisms in the formulas, writing, e.g, $V
  \otimes  U \otimes W$ rather  than $(V \otimes U) \otimes W$. Thus all
  formulas and identities only  make sense after the insertion of
  appropriate canonical morphisms.  Pedantic readers may complete all
  computations by inserting  appropriate canonical morphisms.
\end{remark}

As in the $G = 1$ case, existence of morphism $\delta : V \to V^{**}$
is equivalent to existence of the twist $\theta_V$. Hence we have the
following lemma.
\begin{lemma}\label{l:theta}
  Let $\CC$ be a $G$-equivariant fusion category. Then one can define a
  collection of functorial morphisms $\theta_V: V \to {}^gV$ for $V \in
  \CC_{g}$ which satisfy the following conditions (here $V \in \CC_g$
  and  $U \in \CC_h$):
  \begin{align*}
  \theta_{\1} &= \id \\
  \theta_{U \otimes V} &= R_{{}^{hg}V,{}^hU}R_{{}^hU,{}^gV}(\theta_U
    \otimes \theta_V) \\
  \theta_{V^*} &= R_{g^{-1}}(\theta_{V}^{*}) \\
  \theta_{{}^hV} &= R_{h}(\theta_V).
  \end{align*}

  Conversely, if we have $\theta$ satisfying the above condition then we
  can recover $\delta$ from this $\theta$, $R$, and the monoidal
  structure.
\end{lemma}
\begin{proof}
  The proof is completely parallel to the $G = 1$ case
  discussed in \cite[Section 2.2]{BK}. Details of $G$--equivarint case
  can be found in \cite{equimf1}. 
\end{proof}
\begin{definition}
  Let $\CC$ be a $G$-equivariant category. An object $W$ is called a
  weak dual of $V$ if
  $$
  \Hom(\1,V \otimes X) = \Hom(W,X).
  $$
  In other words,  the functor $\Hom(\1, V \otimes - )$ is represented
  by  the object $W$. In this case, we will  denote $W$ by $V^*$.
  Obviously, in the case of a $G$ equivariant fusion category the usual
  right dual of an object is also a weak dual.
\end{definition}

\begin{definition}
  A $G$-equivariant category $\CC$ is called a $G$-equivariant weakly
fusion category if it satisfies all the conditions
  of  a $G$-equivariant fusion category except for the rigidity
  condition.  Instead of the condition that each object has a right dual
   we require  that each object has a weak dual.
\end{definition}

\begin{remark}
  Of course, a $G$-equivariant fusion category is also a $G$-equivariant
  weakly fusion category but the converse is not true.
\end{remark}

\section{$G$-covers}

The main topological object of our study is the notion of a $G$-cover
of a surface. Detailed description of them is given in the  \cite{TA}.
For readers convenience, we recall basic definitions here.

\begin{definition}\label{d:exsurf}
  An {\em extended surface} is a compact, smooth, oriented, closed
  surface  $\Sigma$ (not necessarily connected), possibly with boundary
  and  with a  choice of  a distinguished (marked) point on
  each of its boundary  components.
\end{definition}
We denote by $A(\Sigma)$ the set of the boundary components of $\Si$. So
an extended surface will be denoted by $(\Sigma,\{p_a\}_{a \in
A(\Sigma)})$ where $p_a$ is the choice of marked point on the $a$ th
boundary component.

\begin{definition}\label{d:Gcover}
  A $G$ cover of $(\Sigma,\{p_a\})$ is a pair 
    $(\pi\colon\tilde{\Sigma} 
  \longrightarrow \Sigma, \{\tpa\})$ where $(\pi\colon 
  \tilde{\Sigma} \longrightarrow \Sigma)$ is a principal $G$-cover
  (possibly not  connected) and $\{\tpa\}$ are choice of points
  in the fiber of $p_a$:  $\tilde{p}_a \in \pi ^{-1}
  (p_a)$  for all  $a \in A(\Sigma)$. 
\end{definition}

For brevity, we will usually denote a $G$-cover just by $\tsi$,
suppressing all other data. 

Note that are $G$-covers are not required to be connected but are
required to be unbranched. 

One can easily define the notion of a morphism between two $G$-covers;
also, since each $G$-cover comes with a choice of a marked points, we
can define, for each boundary component $a$ of $\Si$, the monodromy
$m_a(\tsi)\in G$ (see, e.g., \cite{TA} for details). 
\begin{lemma}\label{l:glue}
  Let  $(\tsi, \{\tilde{p}_a\})$ be a $G$-cover, and let $a,b\in
  A(\Si)$ be two boundary components. Let $\ph\colon (\del \Si)_a\to
  (\del \Si)_b$ be an orientation reversing homeomorphism of boundary
  circles such that $\ph(p_a)=p_b$ \textup{(}it is well-known that such
  a homeomorphism is unique up to isotopy\textup{)}. Then $\ph$ can be
  lifted  to
  an orientation-reversing homeomorphism of boundary components of the
  $G$-cover $\tilde\ph\colon  (\del \tsi)_a\to
  (\del \tsi)_b$ such that $\tilde\ph(\tilde{p}_a)=\tilde{p}_b$ iff the
  monodromy satisfies the following condition:
  \begin{equation}\label{e:glue}
    m_am_b=1.
  \end{equation}
  If this condition is satisfied, then one can form a new $G$-cover by
  identifying  \textup{(}``gluing''\textup{)} $ (\del \tsi)_a\to
  (\del \tsi)_b$ using $\tilde \ph$. We will denote this new $G$-cover
  by 
  $$
    \sqcup_{a,b}\tsi. 
  $$
\end{lemma}  
As a special case, we can consider the situation where $\Si$ is
disconnected: $\Si=\Si_1\sqcup\Si_2$, and $a\in A(\Si_1)$,
$b\in (\Si_2)$. In this case we will use the notation 
$$
\Si_1\sqcup_{a,b}\Si_2
$$
for the result of identification, or ``gluing''.

\begin{definition}
  For every $n\geq 0$, we define the standard sphere, $S_{n}$, to be the
  Riemann sphere $\overline{\mathbb{C}}$ with $n$ disks
  $|z-k|<\frac{1}{3}$ removed and with the marked points being
  $k-\frac{i}{3}$, here $k = 1,2,3,\dots,n$. Of course, we could replace
  these $n$ disks with any other $n$ non-overlapping disks with centers
  on the real line and with marked points in the lower half plane. Any
  two such spheres are homeomorphic and the homeomorphism can be chosen
  canonically up to homotopy. Note that the set of boundary components
  of the standard sphere is naturally indexed by numbers $1,2,\dots,n$.
\end{definition}

\begin{figure}[ht]
  \includegraphics[scale=.6]{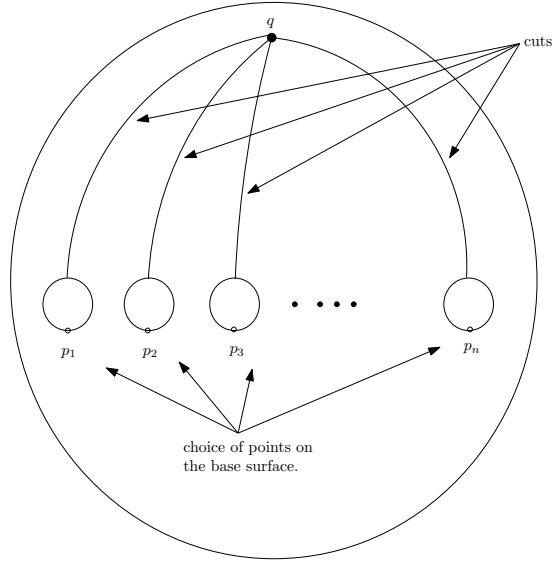}
  \caption{Cuts on $S_{n}$}\label{cutsons11}
\end{figure}

 We will now define the notion of standard blocks, or some
``distinguished'' $G$-covers of the standard sphere. Let us start with
the standard sphere with $n$ holes, $S_{n}$ and $2n$ elements
$\{g_{1},\dots,g_{n}\}$ and $\{h_{1},\dots,h_{n}\}$ in $G$ such that 
 $g_{1}\dots g_{n} = 1$. Make the cuts on $S_{n}$ as
in \firef{cutsons11}.
Here the point $q \in S_{n}$ in \firef{cutsons11} is the point at
$i\infty$. In fact $q$ can be chosen to be any point on the upper
hemisphere as long as it does not belong to the boundary circles. Then
one can easily sees that $S_n \setminus$cuts is simply connected.

Consider the trivial $G$-cover $(S_n \setminus\text{cuts})\times G\to
S_n \setminus\text{cuts}$. Now, define a $G$-cover of $S_n$ by gluing
along $i$-th cut, identifying point $(z,x)\in (S_n
\setminus\text{cuts})\times G$ on the left hand side of the $i$-th cut
with the point $(z,xg_i)$ on the right-hand side. 
One easily sees that this agrees with the
action of $G$ on the fibers (recall that $G$ acts by left
multiplication); condition $g_1\dots g_n=1$ ensures that this gluing
defines a cover which is unbranched at point $q$.
Finally, define a marked point $\tilde{p}_i$ on $i$-th boundray circle
to be $\tilde{p}_i = (p_i,h_i)$. This defines a $G$-cover of $S_n$.

\begin{definition}\label{d:stblock}
  The $G$ cover of $S_{n}$ constructed above will be  called the
  standard block  and will be denoted by $S_n(g_{1},g_{2},\dots,g_{n} ;
  h_{1},h_{2},\dots,h_{n})$. Note that it is only defined if  
  $g_{1}g_{2},\dots,g_{n} = 1$. 
\end{definition}

\begin{lemma}\label{isolemma}
  Let $S_n(g_1,\dots,g_n;h_1,\dots,h_n)$ and
  $S_n(g_1',\dots,g_n';h_1',\dots,h_n')$
  be two standard blocks. Then the identity isomorphism $S_n\to S_n$
  can be lifted to an isomorphism $S_{n}(\mathbf{g,h})
  \to S_{n}(\mathbf{g',h'})$ iff there
  exists $x \in G$  so that $xg_ix^{-1}=g_i'$ and
  $h_ix^{-1}=h_i'$ for $i=1\dots n$; in this  case, the isomorphism is
  unique. We denote the  isomorphism $S_{n}(\mathbf{g,h})
  \to S_{n}(\mathbf{g',h'})$ by  $\phi_{x}$.
\end{lemma}

\begin{proof}
See the paper \cite{TA}.
\end{proof}

\begin{lemma}
  For the standard block  $S_n(g_1,\dots, g_n;h_1,\dots, h_n)$, the 
  monodromy  $m_i \in G$ around the $i$-th boundary circle is given by  
  $$
    m_i =  h_ig_i^{-1}h_i^{-1}.
  $$
\end{lemma}

\begin{proof}
See the paper \cite{TA}.
\end{proof}

\begin{corollary}
  $i$-th boundary circle of $S_{n}(g_{1},...,g_{n};h_{1},...,h_{n})$ can
  be glued to the $j$-th boundary circle of
  $S_{m}(u_{1},\dots,u_{m};v_{1},\dots,v_{m})$ iff
  $h_{i}g_{i}^{-1}h_{i}^{-1} = [v_{j}u_{j}^{-1}v_{j}^{-1}]^{-1}$.
\end{corollary}

\begin{definition}\label{d:param}
  Let $\tsi \to \Sigma$ be a $G$ cover. A {\em parameterization} of
  $\tsi$ an isomorphism
  of   this $G$-cover with one or gluing of several standard blocks:

  $$f\colon
  \tsi \isoto 
  S_{n_{1}}(\bg^1,\bh^1)\sqcup_{i_1,j_1}
  S_{n_{2}}(\bg^2,\bh^2)\sqcup_{i_2,j_2}\dots 
  S_{n_{k}}(\bg^k,\bh^k)
$$`
\end{definition}
It is easy to see that parameterization can be equivalenlty described
by the following data:
 \begin{enumerate}
  \item A finite set $C=\{c_1,\dots \}$  of closed non-intersecting
    curves (``cuts'')   $c_i\in \Si$. 
  \item A choice of marked points $p_c\in c$, one points for each cut
    $c$, and a choice of lifting $\tilde p_c\in \pi^{-1}(p_c)$.
  \item For every connected component $\Si_k$ of $\Si\setminus\{c_i\}$,
    an isomorphism of $G$-covers $f_k\colon \tsi_k\to S_{n_k}(\bg^k,
    \bh^k)$.  
 \end{enumerate}

Finally, note that if $\Si$ is an extended surface with the set of
boundary components $A=A(\Si)$, then we have a natural action of the
group $G^A$ on the category of $G$-covers of $\Si$ by changing the
marked points $\tpa$: if $\bx=\{x_a\}_{a\in A}\in G^A$, then we
define 

\begin{equation}\label{e:tx}
 \bx (\tsi,\{\tpa\})=(\tsi, \{\tilde{p}'_a\}),\qquad 
    \tilde{p}'_a=x_a \tpa.
\end{equation}

\section{$G$-equivariant modular functor}
In this section we introduce the $G$-equivariant analog of the notion
of modular functor; we will call such an analog a {\em $G$-equivariant
modular functor}, or simply  $G$-MF.

Let $\CC$ be a semisimple abelian $G$ equvariant category; 
we assume that the set $I$ of equivalence classes of
simple objects is finite.

\begin{definition}\label{d:R}
Let $\CC$ be a $G$-equivariant category. An object $\RR\in
\CC\tbox\CC$ is called symmetric and $G$-invariant if
\begin{enumerate}
\item 
  $\RR\in \bigoplus_h\CC_h\tbox\CC_{h^{-1}}$
\item $\RR$ is symmetric, i.e. we have an isomorphism  $\si
  \colon \RR\isoto R^{op}$ as in
  \cite[Section~2.4]{BK}. 
\item For every $G$ we have an
  isomorphism $\RR\simeq (R_g\tbox R_g)(\RR)$; these isomorphisms
  should be compatible with each other and with the symmetry isomorphism
  $\si$. 
\end{enumerate}
\end{definition}

As we will show later, a typical example of such an object is when
$\CC$ is a $G$-fusion category and $\RR=\bigoplus V_i\tbox V_i^*$,
where $V_i$ are simple objects. 

We will frequently use the following standard convention: if
$\RR=\bigoplus_i \RR^1_i\tbox \RR^2_i$, then we will drop the index $i$
and summation from our formulas, writing, for example 
\begin{equation}\label{e:Rconvention}
\Hom(A, \RR^1)\otimes \Hom(B,\RR^2)
\end{equation}
for $\bigoplus_i\Hom(A, \RR^1_i)\otimes \Hom(B,\RR^2_i)$.

The following definition is the main definition of this paper; it
generalizes the well-known definition of the modular functor to
$G$-equivariant case.

\begin{definition}\label{d:mf}
  Let $\CC, \RR$ be as above.  A $\CC$-extended $G$-equivariant
  modular functor ($G$-MF for short) is the following collection of
  data:
  \begin{enumerate}
  \item To every $G$-cover $(\tsi, \{\tilde{p}_a\})$ is assigned a
    polylinear functor
    $$
    \tau(\tsi)\colon  \tbox_{a\in A(\Si)}\CC_{m_a^{-1}(\tsi)}\to  \Vec.
    $$
    (Here $m_a$ is monodromy around $a$-th boundary component of the
    $G$-cover).

    In other words, for every choice of objects $W_a\in
    \CC_{m_a^{-1}(\tsi)}$ attached to every boundary component of $\Si$
    is   assigned a finite-dimensional vector space
    $\tau(\tsi;\{W_a\})$, and  this assignment is functorial in $W_a$.

  \item \label{d:mf-2}
    To every morphism of $G$-covers $f\colon \tsi\isoto\tsi'$ is
    assigned a functorial isomorphism $f_*\colon \tau(\tsi)\isoto
    \tau(\tsi')$, which depends only on the isotopy class of $f$.

  \item\label{d:mf-3}
    Functorial isomorphisms $\tau(\emptyset) \isoto k$,
    $\tau(N_1\sqcup N_2) \isoto \tau(N_1)\ttt\tau(N_2)$.

  \item\label{d:mf-4} {\bf Gluing isomorphism:} Let $(\tsi, \{p_a\})$ be
    a $G$-cover and $\al,\be\in A(\Si), \al\ne \be$ be two boundary
    components such that condition \eqref{e:glue} holds. Let
    $\sqcup_{\al,\be}(\tsi)$ be the surface obtained by gluing
    components $\al,\be$ of $\tsi$ as in \leref{l:glue}. Then we have a
    functorial isomorphism
    \begin{equation}\label{e:glue2}
      G_{\al,\be} \colon \tau(\tsi; \{W_a\}, \RR_{\al,\be}) 
      \isoto
      \tau(\sqcup_{\al,\be}\tsi; \{W_a\} ),
    \end{equation}
    where $\RR_{\al,\be}$ means that we assign the symmetric object
    $\RR\in \CC^{\tbox 2}$ to boundary components $\al,\be$.

  \item\label{d:mf-5} For any $G$-cover $(\tsi,\tpa)$ and any
    $\bx=\{x_a\}_{a\in A(\Si)}\in G^{A(\Si)}$,  we have
    functorial isomorphisms
    \begin{equation}\label{e:Tx}
      T_{\bx}\colon 
      \tau(\tsi,\{\tpa\},\{W_a\})\simeq
      \tau(\tsi,\{x_a\tpa\},\{{}^{x_a}W_a\}) 
    \end{equation}

\end{enumerate}

The above data have to satisfy the following axioms:
\begin{description}
\item[Multiplicativity] $(fg)_*=f_*g_*$, $\id_*=\id$.

\item[Functoriality] All isomorphisms in parts \ref{d:mf-3},
  \ref{d:mf-4}, \ref{d:mf-5} above are functorial in $\tsi$.

\item[Compatibility] All isomorphisms in parts \ref{d:mf-3},
  \ref{d:mf-4}, \ref{d:mf-5} above are
compatible with each other.

\item[Symmetry of gluing] After the identification $\RR\simeq
  \RR^{op}$, we have $G_{\al,\be} = G_{\be,\al}$.

\item[Normalization] $\tau(S^2\times G)=k$.
\end{description}
\end{definition}

 Explicit statements of all functoriality and compatibility axioms are
 similar to the ones in $G=\{1\}$ case which can be found
 in \cite{BK}. The only new compatibility relations are those involving
$T_{\mathbf x}$. For the most part, they are quite obvious;  the only
one which is not immediately obvious is the one involving compatibility
of $T_{\mathbf x}$ and gluing, which is given below.

 Let $\tsi$ be a $G$-cover  and $\al,\be \in A(\Sigma_1), \al\ne
\be$. Assume that condition \eqref{e:glue} is 
 satisfied so that we can glue  $\al$ and $\beta$
boundary components. Let $W_a, a\in A'=A(\Si)\setminus\{\al,\be\}$, be a
collection of objects from $\CC$ and let $\bx \in
G^A$ be such that $x_\al=x_\be$.  Then the following diagram commutes:
 $$
 \xymatrix{
 \tau[\tsi, \{\tpa\}; \{W_a\}_{a\in A'},  \RR^1, \RR^2]
   \ar[r]^{G_{\al, \be}} \ar[d]_{T_{\bx}} 
        & \tau[ \sqcup_{\al,\be } \tsi,\{\tpa\}; \{W_a\}]
                           \ar[d]_{T_{\bx '}}\\
 \tau[\tsi, \{x_a\tpa\}; \{{}^{x_a}W_a\}, {}^x\RR^1, {}^x\RR^2]
\ar[r]_{G_{\al,\beta}}& 
    \tau[ \sqcup_{\al,\be } \tsi,\{x_a\tpa\}; \{{}^{x_a}W_a\}]
}
$$

 Here we assigned $\RR^1$ and $\RR^2$ to $\al$ and $\be$ boundary
components respectively; the bottom arrow also uses isomorphism
${}^x\RR^1\tbox {}^x\RR^2\simeq \RR^1\tbox \RR^2$ (see \deref{d:R}). 

Also note that condition $x_\al=x_\be$ ensures that if we can glue the
$\al$ and $\be$ boundary components of $(\tsi, \{\tpa\})$ then we can
also glue the $\al$ and $\be$ boundary components of
$(\tsi,\{x_a\tpa\}) $.

 We leave all the other compatibility conditions, which are rather 
obvious, to the imagination of the readers.

\begin{definition}\label{d:nondegenerate}
 A $G$-MF is called non-degenerate if for every non-zero object $X\in
  \CC$, there exists a $G$-cover $(\tsi, \{\tpa\})$ and a collection of
objects 
  $\{W_i\}$ such that $\tau(\tsi, \{\tpa\}; \{X, W_1, \dots \})\ne 0$. 
\end{definition}
As in $G=\{1\}$  case, it is sometimes convenient to consider more
restricted version of modular functor, in which we only allow genus 0
surfaces.

\begin{definition} \label{d:mf0}
A genus zero $\CC$ extended $G$-equivariant modular functor
  (genus 0 $G$-MF for short) consists of the same
  data as defined in \deref{d:mf} except that $\tau$ is only defined
  for surfaces all connected components of which have genus zero, and
  the gluing is only defined if $\al,\be$ are in different
  connected components of $\Si$.
\end{definition}

In this paper, we will only consider genus zero $\CC$ extended
$G$-eqivariant modular functor; case of positive genus will be
addressed in forthcoming papers. 

\section{Statement of the main theorem}
The following is the main theorem of this paper. The rest of the paper
is
devoted to prove this theorem.
\begin{theorem*}
  Let $\CC$ be a semisimple $G$-equivariant abelian
  category $\CC$ and let 
  $$
  \{V_{i} \st   i \in I \}
  $$  
  be a set of representative of isomorphism classes of simple objects.
  Assume that $|I| < \infty$.

  \begin{enumerate}
    \item If we have a non-degenerate $\CC$-extended  genus zero $G$-MF,
      then we can construct a structure of $G$-equivariant weakly
      fusion category on $\CC$. 
    \item Conversely, if we have a structure of $G$ equivariant fusion
      (or weakly fusion) category on $\CC$ then we can define a
      non-degenerate $\CC$ extended genus zero $G$-MF. 
  \item The above two constructions are inverse of each other. 
  \end{enumerate}
  The above correspondence between modular functor and fusion category
  satisfies the following properties:
  \begin{enumerate}
  \item The $G$-invariant symmetric object $\RR$ used in the definition
      of the MF is given by $\RR=\bigoplus V_i\tbox V_i^*$.
  \item Let $S_{n}(\mathbf{g,h})$ be the standard block as defined in
    \deref{d:stblock}. Then 
  $$
  \tau[S_{n}(\mathbf{g,h});W_1,\dots, W_n)] =
  \Hom_{\CC}(\1,{}^{h_1^{-1}}W_1 \otimes\dots \otimes {}^{h_n^{-1}}W_n).
  $$

  \item Let  $\tilde{z}$ be the isomorphism between
    standard    blocks given by the rotation (see 
    \cite{TA} for details):
    \begin{equation}\label{e:z}
      \tilde{z}\colon S_n(g_1,\dots,g_n;h_1,\dots,h_n) \to
      S_n(g_n,g_1,\dots, g_{n-1};h_n,h_1,\dots,h_{n-1}).
    \end{equation}
    Then $\tilde{z}_*$ is given by:
   \begin{align*}
    &\Hom_{\CC}(\1, X_1 \otimes \dots\otimes X_n) \to
    \Hom_{\CC}(\1, X_1 \otimes \dots\otimes X_n^{**}) \\
    &\qquad \to \Hom_{\CC}(X_n^*,X_1\otimes\dots\otimes X_{n-1})
    \to \Hom_{\CC}(\1, X_n \otimes X_1\otimes\dots\otimes X_{n-1})
    \end{align*}
    where $X_i={}^{h_i^{-1}}W_i$. 
  \item Let $\tilde{b}$ be the braiding morphism between
    standard  blocks \textup{(}see \cite{TA} for details\textup{)}:
    $$
    \tilde{b}\colon S_3(g_1,g_2,g_3;h_1,h_2,h_3) \to
      S_3(g_1,g_2g_3g_2^{-1},g_2;h_1,h_3g_2^{-1},h_2)
    $$
    Then $\tilde{b}_*$ is given by
    \begin{align*}
    \Hom_{\CC}&(\1,X_1\otimes X_2\otimes X_3) 
       \xxto{\id \otimes R_{X_2,X_3}} 
        \Hom_{\CC}(\1,X_1\otimes {}^{g_2}\!X_3\otimes X_2)\\
    &=\Hom_{\CC}(\1,\otimes {}^{h_1^{-1}}W_1\otimes {}^{g_2h_3^{-1}}W_3
      \otimes {}^{h_2^{-1}}W_2 )
    \end{align*}
    where, as before, $X_i={}^{h_i^{-1}}W_i$. 
    Note that since $W_2 \in \CC_{h_2g_2h_2^{-1}}$,
    $X_2={}^{h_2^{-1}}W_2  \in \CC_{g_2}$. So by the definition of
    the braiding isomorphism of  the fusion   category, we need to
    twist $X_3$ by $g_2$.

  \item Let 
    $$
    \phi_x\colon  S_n(g_1,...,g_n;h_1,...,h_n) \to
    S_n(xg_1x^{-1},...,xg_nx^{-1};h_1x^{-1},...,h_nx^{-1})
    $$
    be the isomorphism between standard blocks described in
    \leref{isolemma}.  Then $(\phi_x)_*$ is
    given by the following formula:
    \begin{align*}
    \Hom_{\CC}& (\1,{}^{h_1^{-1}}V_1 \otimes\dots\otimes
                    {}^{h_n^{-1}}V_n)     \to
    \Hom_{\CC}({}^x\1,{}^{xh_1^{-1}}V_1 \otimes \dots\otimes
      {}^{xh_n^{-1}}V_n)\\
    &\to
    \Hom_{\CC}(\1,{}^{xh_1^{-1}}V_1 \otimes \dots\otimes
    {}^{xh_n^{-1}}V_n)
    \end{align*}
    We use the fact that ${}^x\1 \cong \1$ and the action of the group
    $G$ is a   tensor functor.
  \end{enumerate}
\end{theorem*}

  This is not the full list of all the properties enjoyed by the
correspondence between modular functors and equivariant categories; more
properties will be seen from the construction. 
\begin{proof}[Idea of the proof]
   Instead of proving a direct correspondence between
   $G$-equivariant modular functor and $G$-equivariant fusion category,
   we will introduce an intermediate object, which we will call  $G$
   equivariant Moore-Seiberg data and then established the
   equivalence between these three notions:
\begin{equation}\label{e:idea}
    \parbox{1in}{$G$-eqivariant\\ modular functor}
       \longleftrightarrow\quad
\parbox{0.9in}{$G$-equivariant
      MS data}
\longleftrightarrow\quad 
\parbox{1.5in}{$G$-equivariant\\ weakly fusion category}
\end{equation}
\end{proof}

\section{$G$-equivariant Moore-Seiberg data}
In this section we wil introduce the intermediate object, the
Moore-Seiberg data  (MS for short); the goal of this  is
encoding the structure of a fusion category in terms of vector spaces 
$\<W_1,\dots, W_n\>=\Hom_\CC (\1, W_1\otimes \dots\otimes W_n)$ and
suitable isomorphisms between such spaces. 

For Moore-Seiberg data for $G = 1$ case, see  \cite[Section 
5.3]{BK}. Here 
we need the extension of this concept to $G$-equivariant case. 

\begin{definition}\label{d:MSdata}
Let $\CC$ be  a $G$-equivariant abelian  category (see 
\deref{geqcat}). Then
$G$-equivariant Moore-Seiberg data ($G-$MS data for short) is the
following collection of data:
\begin{description}

  \item[Conformal blocks]\label{ms1}
     For each $n \geq 0$ and $m_1,m_2,\dots,m_n
      \in G$ satisfying  $m_1m_2\cdots m_n = 1$  we have a functor:
      $$
      \< {} \>\colon \CC_{m_1} \tbox\cdots\tbox \CC_{m_n}
      \rightarrow \Vec
      $$
    where  $\Vec$ denotes the category of finite-dimensional vector
    spaces. 

    \textbf{Note}: We can trivially extend this functor to a functor
    $\CC\tbox \cdots \tbox \CC\to \Vec$ by letting
      $\< W_1,W_2,\dots, W_n \> = 0$ if $W_i\in m_i$,  $m_1m_2...m_n
      \neq 1$.

  \item[$G$-invariance]\label{ms2}
    For each $g \in G$, we have functorial isomorphism
    $$
    \phi_{g}\colon \< V_1,V_2,\dots, V_n \> \longrightarrow
    \< {}^gV_1,{}^gV_2,\dots, {}^gV_n \>.
    $$
    satisfying $\phi_{gh}=\phi_h\phi_g, \phi_1=\id$. 
    
    \item[Rotation isomorphism]\label{ms3}
      Functorial isomorphism
      $$
        Z\colon \< V_1,V_2,\dots, V_n \> \rightarrow \<
        V_n,V_1,\dots,V_{n-1} \>
      $$
    \item[Symmetric object]\label{ms4}
      A symmetric $G$-invariant object $\RR\in \CC\tbox\CC$ as in
      \deref{d:R}.

    \item[Gluing isomorphism]\label{ms5}
     For each $k,l \in \Z_{+}$,
    there exist functorial isomorphism
    $$
    \G\colon  \< A_1,\dots,A_k,\RR^1 \> \otimes \<
      \RR^2,B_1,\dots,B_l
    \> \rightarrow \< A_1,\dots,A_k,B_1,\dots,B_l \>
    $$
    (As before, $\RR^1, \RR^2$ should be understood as in
    \eqref{e:Rconvention}). 

    \item[Commutativity isomorphism]\label{ms6}
      For $A \in \CC_g$ and $B \in \CC_h$, we have functorial
      isomorphism
        $$
        \sigma\colon \<X,A,B \> \to \<X, {}^{g}\!B,A\>. 
        $$
        \textbf{Note}: if $X \in \CC_{p}$ then $\<X,A,B \>$
          is zero unless $pgh = 1$. But then ${}^{g}\!B \in
          \CC_{ghg^{-1}}$
          and the  right hand side is zero unless we have 
          $pghg^{-1}g =  pgh = 1$; thus, if one side is zero for
          grading reasons, then so is the other.
      \end{description}

      These above data must satisfy the axioms formulated below:
non-degeneracy, normalization, rotation axiom,  associativity of $G$,
symmetry of $G$, compatibility of $\phi$, hexagon axiom, and Dehn twist
axiom. 

\end{definition}

Before formulating the axioms, it would be convenient to define certain
compositions of elementary gluing, rotation, and commutativity
isomorphisms as follows. 

\begin{description}
  \item[Generalized gluing]
    For any $k,l,m\ge 0$, we define {\em generalized gluing}
    isomorphism
    \begin{equation}\label{e:gen_gluing}
    \begin{aligned}
      \G&\colon 
        \<A_1,\dots, A_k, \RR^1, C_1,\dots, C_m\> 
           \otimes \< \RR^2, B_1,\dots,B_l\>\\
      &\quad \to \<A_1,\dots, A_k, B_1,\dots, B_l, C_1,\dots, C_m\>
    \end{aligned}    
    \end{equation}
    as the following compoisition:
    \begin{align*}
      &\<A_1,\dots, A_k, \RR^1, C_1,\dots, C_m\> 
           \otimes \< \RR^2, B_1,\dots,B_l\>\\
         &\qquad   \xxto{Z^{m}\otimes \id}
       \<C_1,\dots, C_m,A_1,\dots, A_k, \RR^1\> 
           \otimes \< \RR^2, B_1,\dots,B_l\>\\ 
       &\qquad    \xxto{\G}
       \<  C_1,\dots, C_m, A_1,\dots, A_k,B_1,\dots, B_l\>\\ 
       &\qquad   \xxto{Z^{-m}}
      \< A_1,\dots, A_k,B_1,\dots, B_l, C_1,\dots, C_m \> 
      \end{align*}

    \item[Generalized commutativity]
       For any $k,l\ge 0$, we define the generalized commutativity
      isomorphisms 
        \begin{equation}\label{e:gen_commut}
        \sigma\colon \< A_1,\dots,A_k,X,Y,B_1,\dots,B_l \>  
          \to \<A_1,\dots,A_k,{}^pY,X,B_1,\dots,B_l \>
        \end{equation}
        where $X \in \CC_p$ and $Y \in \CC_q$, as the following
        composition 
        \begin{align*}
         &\< A_1,\dots,A_k,X,Y,B_1,\dots,B_l \> \xxto{\G^{-1}}
          \< A_1,\dots,A_k,\RR^1,B_1,\dots,B_l \>
           \otimes \<\RR^2, X,Y\>\\
         &\qquad\xxto{\id\otimes \si}          
             \<A_1,\dots,A_k,\RR^1,B_1,\dots,B_l\>
              \otimes \<\RR^2,{}^pY,X\>\\
          &\qquad\xxto{\G}
          \<A_1,\dots,A_k,{}^pY,X,B_1,\dots,B_l \>
        \end{align*}
        (here $\G$ is the generalized gluing \eqref{e:gen_gluing}).
         
        From now on we will denote by $\sigma$ both the generalized
        commutativity isomorphisms and the usual commutativity
        isomorphisms.

\end{description}
  One can define even more general isomorphisms; however, it won't be
  necessary for our purposes.

Now we are ready to formulate the axioms of the MS data.

    \begin{description}
      \item[Non-degeneracy]\label{msax1}
        For each  object $X \in \CC$, there exists an object $X
          \in \CC$ so that $\< X,V \>\neq 0$.
      
        \item[Normalization]\label{msax2}{For $n =0$
          $$
          \< \>\colon \Vec \to \Vec
          $$
          is the identity functor (as before, $\Vec$ denotes  the
          category of finite-di\-men\-si\-on\-al vector spaces).
    
        \item[Associativity of $\G$]\label{msax3}
         Let $\RR, \tilde \RR$ be two copies of $\RR$. Then the diagram
          in \firef{associativity_of_g} is commutative. 
        \begin{figure}[ht]
         $$\xymatrix@C=-80pt{
              & \scriptstyle{\<A_1,\dots,A_k, \RR^1\>
                      \otimes\<\RR^2,B_1,\dots, B_l, \tilde \RR^2\>
                      \otimes\<\tilde \RR^2,C_1,\dots, C_p\>}
                      \ar[dl]_{\G\otimes\id} \ar[dr]^{\id\otimes \G}\\
               \scriptstyle{
               \<A_1,\dots,A_k, B_1,\dots, B_l, \tilde \RR^2\>
                      \otimes\<\tilde \RR^2,C_1,\dots, C_p\>}
                      \ar[dr]_{\G} &&
               \scriptstyle{\<A_1,\dots,A_k, \RR^1\>
                      \otimes\<\RR^2,B_1,\dots, B_l, C_1,\dots, C_p\>}
                      \ar[dl]^{\G} \\
                &\scriptstyle{
                \<A_1,\dots,A_k, B_1,\dots, B_l,C_1,\dots, C_p\>}
          }
          $$
          \caption{Associativity axiom for $\G$}
          \label{associativity_of_g}
        \end{figure}

      \item[Rotation axiom]\label{msax4}
        The isomorphism of the vector spaces
        $$
          Z^n\colon \< V_1,\dots,V_n \> \to 
                \< V_1,\dots,V_n \>
        $$
        is equal to  identity.
    
      \item[Symmetry of $\G$]\label{msax5}
         Again let $\RR = \RR^1 \tbox \RR^2$
        and $P$ be the usual isomorphism between 
        the vector spaces $A\otimes B \to
          B \otimes A$ given by $P(a \otimes b) = b \otimes a$. Then
          the diagram in \firef{ms_g_symmetry} is commutative. 

        \begin{figure}[ht]
        $$\xymatrix{
           \<A_1,\dots,A_k, \RR^1\>\otimes 
                    \<\RR^2,B_1,\dots, B_l\> \ar[r]^\G 
                        \ar[d]_{Z\otimes Z^{-1}}
               & \<A_1,\dots,A_k, B_1,\dots, B_l\>\ar[dd]^{Z^l}\\
          \<\RR^1, A_1,\dots,A_k \>\otimes 
                    \<B_1,\dots, B_l, \RR^2\>\ar[d] \\
           \<B_1,\dots, B_l, \RR^2\>\otimes
                     \<\RR^1, A_1,\dots,A_k\>\ar[r]^\G
           & \<B_1,\dots, B_l,  A_1,\dots,A_k\>
          }
          $$
          \caption{Symmetry of $\G$ axiom}\label{ms_g_symmetry}
        \end{figure}

      \item[Compatibility of $\phi$]\label{msax7}
        The functorial isomorphisms $\phi_g$ must be compatible with all
        the other isomorphisms. More precisely we
        have the following:
        \begin{enumerate}
        \item 
           $\phi$ must be compatible with rotation isomorphism, i.e.   
          the following diagram is commutative:
            $$
            \xymatrix{
             \<W_1,\dots, W_n\>\ar[r]^{Z} \ar[d]^{\phi_g} &
                  \<W_n, W_1,\dots, W_{n-1}\>\ar[d]^{\phi_g}\\
             \<{}^gW_1,\dots, {}^gW_n\>\ar[r]^{Z}  &
                  \<{}^gW_n, {}^gW_1,\dots, {}^gW_{n-1}\>
            }
           $$
        \item 
          $\phi$ must be compatible with the gluing isomorphism, i.e.
            for  each $k,l \in \Z_+$ the following diagram  must be
            commutative
            $$
            \xymatrix{
            \<A_1,\dots, A_k,\RR^1\>\otimes 
                \<\RR^2, B_1,\dots, B_l\>
                \ar[r]^{\G}\ar[d]^{\phi_g\otimes\phi_g} &
            \<A_1,\dots, A_k, B_1,\dots, B_l\> 
                \ar[dd]^{\phi_g} \\
            \<{}^gA_1,\dots, {}^gA_k,{}^g\RR^1\>\otimes 
                \<{}^g\RR^2,{}^gB_1,\dots, {}^gB_l\>
                \ar[d] \\
            \<{}^gA_1,\dots, {}^gA_k,\RR^1\>\otimes 
                \<\RR^2,{}^gB_1,\dots, {}^gB_l\>
                \ar[r]^{\G} &
             \<{}^gA_1,\dots, {}^gA_k, {}^gB_1,\dots, {}^gB_l\>
            }
            $$
                The diagram above uses the isomorphism 
            ${}^g\RR^1 \tbox  {}^g\RR^2 = \RR^1 \tbox \RR^2$.

        \item 
          $\phi$ must be compatible with commutativity isomorphism,
          i.e. for  any  $A \in \CC_p$ and $B \in \CC_q$
          the following diagram must be commutative.
           $$\xymatrix{
             \<X, A, B\>\ar[r]^\si\ar[d]_{\phi_g}
                & \<X, {}^pB, A\>\ar[d]_{\phi_g}\\
             \<{}^gX, {}^gA, {}^gB\>\ar[r]^\si
                & \<{}^gX, {}^{gp}B, {}^gA\>
           }
           $$

        \end{enumerate}}

    \item[Hexagon axiom]\label{msax8} 
      For any $A \in \CC_p, B \in \CC_q$ and $C
        \in \CC_r$, the  diagram below  and a similar diagram with
        $\si$ repalced by $\si^{-1}$ must be 
      commutative 
      $$
       \xymatrix{
        \<X, A,B,C\> \ar[rr]^{\si_{A,BC}} \ar[dr]_{\si_{A,B}}
               &&\<X, {}^pB,{}^p C,A\>\\
        & \<X, {}^pB,A, C\>\ar[ur]_{\si_{A,C}}
        }
       $$     
     
     where $\si_{A,B}$, $\si_{A,C}$ are generaized braidings defined
      by \eqref{e:gen_commut} and $\si_{A,BC}$ is defined as the
      following composition
    \begin{align*}
    &\<X, A,B,C\>\to \<X, A, \RR^1\>\otimes \<\RR^2,B,C\>\\
    &\to \<X, {}^p\RR^1,A\>\otimes \<\RR^2,B,C\>
      \to \<X, \RR^1,A\>\otimes \<{}^{p^{-1}}\RR^2,B,C\>\\
    &\to \<X, \RR^1,A\>\otimes \<\RR^2,{}^pB,{}^pC\>
      \to \<X, {}^pB,{}^p C,A\>
    \end{align*}

    \item[Dehn twist axiom]\label{msax9}
      Let $A \in \CC_p$ and $B \in \CC_q$ with $pq = 1$. Then the Dehn
      twist axiom is the commutativity of the
      diagram in \firef{f:ms_dehn_twist}. Note that if $pq \neq 1$
        then all the vector spaces in \firef{f:ms_dehn_twist} will be
        zero and the diagram  trivially commutes.

      \begin{figure}[ht]
      $$ 
       \xymatrix@R=5pt{
       \<A, B\> \ar[r]^\si \ar[dd]_Z & 
                       \<{}^pB,A\>\ar[dr]^Z\\
        & & \<A, {}^p B\>\\
       \<B,A\>\ar[r]^\si &\<{}^qA, B\>\ar[ur]_\ph                       
       }$$

     \caption{Dehn twist axiom for $G$-MS data.}\label{f:ms_dehn_twist}
      \end{figure}

\end{description}


\section{$G$-MS data and $G$-fusion categories}
We will divide the proof of the main theorem into two steps, first
relating the notion of  $G$-MS data with $G$-equivariant fusion
categories, and then relating
$G$-MS data with the $G$-modular functor. In this section we do the
first step, showing  the equivalence of $G$-equivariant (weakly) fusion
category and $G$-equivariant MS data. 
\begin{theorem}\label{fums}
  \begin{enumerate}
    \item If we have a structure of $G$ equivariant (weakly) fusion
      category on $\CC$ then from this we can create a $G$ equivariant
      MS data on $\CC$.
    \item Conversely, given a $G$ equivariant MS data on $\CC$, we can
      define a structure of $G$ equivariant weakly fusion category on
      $\CC$.
    \item The above two constructions are inverse of each other.
  \end{enumerate}
\end{theorem}
The proof of this theorem is given below. For the most part, it is
parallel to the proof in $G=\{1\}$ case, given in \cite[Section
5.3]{BK}; thus we will only provide detailed explanations of the steps
which are new to $G$-equivariant case. 

\subsection{From fusion categories to $G$-MS data}
Assume that $\CC$ has a structure of $G$ equivariant
(weakly) fusion category. We define the $G$ equivariant MS data as
follows
\begin{description}
  \item[Conformal blocks]
      Let  $W_1 \in \CC_{m_1}, \dots,W_n \in
        \CC_{m_n}$, where $m_1m_2\dots m_n = 1$.  Then define
      \begin{equation}\label{e:conf_blocks2}
    \< W_1,W_2,\dots,W_n \> = \Hom(\1,W_1  \otimes \dots\otimes W_n)
      \end{equation}

  \item[$\phi$-axiom] 
    For $g \in G$, define 
    $$
    \phi_g\colon  \< W_1,W_2,\dots,W_n \> \to 
     \< {}^{g}W_1,{}^{g}W_2,\dots,{}^{g}W_n \>
    $$
    by
    $$
    \Hom(\1,W_1  \otimes \dots,  \otimes W_n)
    \xxto{R_g} \Hom({}^g\1,{}^g(W_1 \otimes \dots \otimes W_n)) 
     \to \Hom(\1,{}^gW_1  \otimes \dots \otimes {}^gW_n)
    $$
    Here $R_g$ is the action of the element $g \in G$. Since $R_g$ is 
    by definition a   tensor functor, we have canonical
    isomorphisms ${}^{g}\1 = \1$ and 
    ${}^{g}(W_1 \otimes \dots \otimes W_n) 
    = {}^{g}W_1 \otimes \dots \otimes {}^{g}W_n$

  \item[Rotation isomorphism]
    Define the rotation isomorphism 
      $$
      Z\colon \< W_1,W_2,\dots,W_n \> \to \<W_n,W_1,\dots,W_{n-1}\>
      $$
    by
    \begin{equation}
    \begin{aligned}
    &\Hom(\1,W_1 \otimes W_2 \otimes \dots \otimes W_n) 
      \to \Hom(\1, W_1\otimes W_2 \otimes \dots\otimes W_{n-1}
                   \otimes W_{n}^{**})\\ 
     &\quad  \to \Hom(W_n^{*},W_1 \otimes\dots\otimes W_{n-1}) 
       \to \Hom(\1,W_n \otimes W_1 \otimes\dots\otimes W_{n-1})
    \end{aligned}
    \end{equation}

    Here we use the  rigidity isomorphisms, balancing isomorphism
    $\delta\colon V \to
    V^{**}$, and the fact that in any weakly rigid fusion category we
    have an isomorphism $\Hom(1, X\otimes W^{**})\simeq \Hom(W^*,X)$
    (see \cite{BK}). 

  \item[Symmetric object]

    The symmetric object $\RR$ is defined by 
    $$
    \RR = V_i \tbox V_{i}^{*}
    $$
    where $\{V_i\}_{i \in I }$ are representatives of the isomorphism
    classes of simple objects. The fact that it is independent
    of the choice of representatives and symmetry of $\RR$ are shown
   in the same  way as in \cite{BK}. 

    To show $G$-invariance, note that since $R_g$ is a tensor functor,
      it must take a simple object to a simple
    object. Thus,  for every $i$ we can choose an
    isomorphism $\psi_g\colon R_g(V_i)\to V_j$ for some $j\in I$. Now
    define isomorphism
    $$
    (R_g\tbox R_g)\RR=\bigoplus R_g(V_i)\tbox R_g(V_i^*)
        \simeq \bigoplus R_g(V_i)\tbox (R_g(V_i))^*
        \xxto{\psi_g\tbox (\psi_g^*)^{-1}} 
        \bigoplus V_j\tbox V_j^*
    $$
    It is easy to see that this isomorphism does not depend on the
    choice of $\psi$ and is compatible with the symmetry $\si\colon
    \RR\to\RR^{op}$.

\item[Gluing isomorphism]
  The gluing isomorphism
  $$
    \< A_1,\dots ,A_k,\RR^1 \> \otimes \< \RR^2,B_1,\dots,B_l \>
    \to \< A_1,\dots ,A_k,B_1,\dots,B_l \>
  $$
  is defined by 
  \begin{equation}\label{e:gluing2}
  \begin{aligned}
  &\bigoplus \Hom(\1,A_1 \otimes\dots\otimes A_k \otimes V_i^*) 
       \otimes \Hom(\1,V_i\otimes B_1 \otimes\dots\otimes B_l)\\
  &\quad \cong
     \bigoplus \Hom(\1,A_1 \otimes\dots\otimes A_k \otimes V_i^*) 
       \otimes \Hom(V_i^*, B_1 \otimes\dots\otimes B_l)\\
  &\quad \cong
     \bigoplus \Hom(\1,A_1 \otimes\dots\otimes A_k\otimes   B_1
   \otimes\dots\otimes B_l)
  \end{aligned}
 \end{equation}

  \item[Commutativity isomorphism] 
   Let $A \in \CC_p$ and $B \in \CC_q$. Define the commutativity
  isomorphism
    $$
    \si\colon \< X, A,B \> \to \<X, {}^p B, A\>
  $$
  by 
  \begin{equation}\label{e:commutativity2}
  \begin{aligned}
   & \Hom(\1, X \otimes A \otimes B) \xxto{\id\otimes R_{A,B}} 
     \Hom(\1,X \otimes {}^pB \otimes A).
  \end{aligned}
  \end{equation}
\end{description}
Then the above data satisfies the definition of $G$-MS data. Indeed, it
can easily be shown that the isomorphisms defined above satisfy all the
axioms of MS-data; the easiest way to do this is to use the technique
of using appropriately marked ribbon graphs to represent morphisms in a
$G$-equivariant fusion category (see \cite{T2}, \cite{orbi1}). 

\subsection{From $G$-MS data to $G$-fusion categories}
Now we are assuming that we have a $G$-equivariant Moore-Seiberg data
and from this we want to construct a structure of $G$ equivariant weakly
fusion category on $\CC$. This construction is parallel to the
construction in the case of $G = \{1\}$ given in \cite[Section 5.3]{BK}.

We will construct the structure of weakly fusion category step by step
as follows:
\begin{description}
\item[Duality]
 Define the functor $*$ by
    $$
    \< V,X \> = \Hom(V^*,X)
    $$
    \begin{remark}
    Every object $V$ of $\CC$ is completely determined by the functor
$\< V, . \>$.
    \end{remark}
   The arguments in \cite[Section 5.3]{BK}, without any changes, show
  that $*$ is an antiequivalence of categories, and that one has a
  canonical isomorphism $\RR\cong \bigoplus V_i\tbox V_i^*$.

\item[Construction of tensor product]
  As in \cite{BK}, define tensor product functor by  
  $$
    \<X, A\otimes B\>=\<X, A, B\>
  $$
  Note that for $X\in \CC_r,A\in \CC_p, B\in \CC_q$, one has $\<X, A,
  B\>=0$ unless $rpq=1$; this implies that $A\otimes B\in \CC_{pq}$. 

  The same arguments as in \cite[Section 5.3]{BK} show that this tensor
  product  has a canonical associativity isomorphism; moreover, we have
  canonical isomorphisms 
    $$
    \<X, A_1 \dotimes A_n \> \cong  \< X, A_1,\dots,A_n \>.
    $$
   
\item[Construction of unit object]

  We define the unit object by 
  $$
      \< \1,X \> = \< X \>
  $$
  Again, the same argument as in \cite[Section 5.3]{BK}, with no changes
  at all, shows that so defined $\1$ is a unit object with respect to
  previously defined tensor product, and that we have a canonical
  isomorphism $\1^*\cong \1$. 

\item[Construction of commutativity isomorphism]

  Let $A \in \CC_p$ and $B\in \CC_q$. We define
  $$
    R_{A,B}\colon  A \otimes B \to {}^pB \otimes A
  $$
  as the following composition:
  $$
    \< X,A \otimes B \> = \< X,A,B \>
    \xxto{\si} \< X,{}^pB,A \> = \< X,{}^pB
    \otimes A \>
  $$
  Here $\sigma$ is the commutativity isomorphism of $G$ equivariant MS
  data.

  Easy explicit computation shows that  the hexagon axiom of fusion
  category is exactly equivalent to the hexagon axiom of $G$ equivariant
  MS data.

\item[Construction of balancing isomorphims]

  We know that having functorial isomorphism
$\delta\colon  V \to V^{**}$ is equivalent to having functorial
isomorpisms (or twist) $\theta_V\colon  V \to {}^gV$ for $V \in \CC_g$
satisfying certain conditions (see \leref{l:theta}). 

To define $\th$, recall the generalized commutativity isomorphisms
defined by  \eqref{e:gen_commut}. In particular, letting $k,l=0$, we
get a commutativity isomorphism 
    $$
    \si\colon \< A,B \> \to  \< {}^pB,A \>, \qquad 
    A \in \CC_p,\ B \in \CC_q. 
    $$
    Now define the  twist functor $\th_V\colon  V \to {}^gV$, $V\in
    \CC_g$ as follows. For any $X\in \CC$, consider the composition 
   
    $$
    \< V,X \>\xxto{\si^{-1}} \< X,{}^gV \> \xxto{Z} \<{}^gV,X \>.
    $$
    This gives the functorial isomorphism between the functor $\< V,.\>
    $ and $\< {}^gV,.\> $ which in turn gives the twist $\th_V$. 
    
    All the required properties of $\th$, listed in \leref{l:theta},
    now easily follow from the properties of commutativity and
    associativity morphisms in teh definition of MS data. 

\end{description}

\section{MS data and  modular functor}
 In this section we do the second step of the proof, showing  the
equivalence of $G$ equivariant genus
zero modular functor and $G$ equivariant MS data. 
\begin{theorem}\label{t:mfms}
  Let $\CC$ be a semisimple $G$-equivariant abelian category with a
  finite number of equivalence classes of simple objects.
  \begin{enumerate}
  \item If we have a non-degenerate $G$ equivariant $\CC$ extended genus
    zero modular functor then we can define a $G$ equivariant MS data
    on  $\CC$
  \item  Conversely, given a $G$ equivariant MS data on $\CC$, we can
    define a non-degenerate $G$ equivariant $\CC$ extended genus zero
    modular functor.
\item The above two constructions are inverse to each other. 
\end{enumerate}
\end{theorem}
The proof of this theorem is given in two subsections below. 

\subsection{From  $G$-MS data to $G$-MF}
Now we are assuming that we have a $G$ equivariant MS data on $\CC$ and
from this we want to construct a $G$ equivariant $\CC$ extended genus
zero modular functor. The construction is similar ot the one given i
n \cite{BK} in $G=\{1\}$ case: we first define the modular
functor on standard blocks. Since any $G$ cover is isomorphic to gluing
of several standard blocks (this identification is not unique; in fact
there are infinitely many parameterization of a given $G$ cover), this
will give us the modular functor on parametrized  $G$-covers. After
this, we show that the modular functor spaces obtained from  any  two
parameterizations of the same $G$ cover are canonically isomosrphic. 

\begin{definition}\label{d:mf_1}
Given $W_{i} \in \CC_{h_ig_ih_i^{-1}}$ for $i=1\dots n$, we define
$$
\tau \left[ S_{n}(g_1,\dots ,g_n;h_1,\dots,h_n; W_1,\dots,W_n)
\right] = \<X_1,\dots,X_n \>
$$
where for brevity we denoted 
$$
X_i={}^{h_{i}^{-1}}W_i.
$$
\end{definition}
\begin{remark}
Note that monodromy  around $i$th boundary component of\\ 
$S_{n}(g_1,\dots,g_n;h_1,\dots,h_n)$ is $m_i=h_ig_i^{-1}h_i^{-1}$,
so we have $W_i\in \CC_{m_i^{-1}}$, as required in the definition of
the modular functor. Note also that $X_i\in \CC_{g_i}$,
so condition $g_1\dots g_n=1$ given by definition of standard block
matches the condition required for $\<X_1,\dots,X_n  \>$ to
be non-zero. 
\end{remark}

Next, let $\tsi$ be a $G$-cover and  $W_a\in \CC_{m_a^{-1}}, a\in
A(\Si)$ --- a collection of objects assigned to boundray components of
$\Si$. Let 
$$
f\colon  \tsi \to S_{n_{1}}(\bg^1,\bh^1)
\bigsqcup_{glued} \dots   \bigsqcup_{glued}
S_{n_k}(\bg^k,\bh^k)
$$
be a parameterization of $\tsi$ (see \deref{d:param}). For each
connected component $\Si_i$ of $\Si\setminus\text{cuts}$, let us assign
an object of $\CC$ to each boundary component of $\Si_i$, by putting
$W_a$ on the boundary of $\Si$ and a copy of $\RR$ on each cut (i.e.,
assigning $\RR^1$ to the component on one side of the cut and $\RR^2$ on
the other side; since $\RR$ is symmetric, choice of the side is not
important).

Then we define the modular functor for parameterized $G$-covers by 
$$
\tau [ \tsi,\{W_a\},f] =
\bigotimes_i \tau[S_{n_{i}}(\bg^i,\bh^i)]
$$
with the choice of objects as explained above. 
This defines a modular functor for parameterized surfaces; by
definition, it satisfies the gluing axiom. Now we need to identify
the spaces $\tau [ \tsi,\{W_a\},f]$ for different parameterizations. 

To do so, we use the same strategy used in $G=\{1\}$ case in
\cite{BK}. Namely, recall  the complex $\M$ from
\cite{TA}; vertices of this complex are exactly different
parameterizations of $\tsi$, edges are certain ``simple moves'' between
parameterizations, and 2-cells describe relations. The main result of
the paper \cite{TA} is that under these basic moves and relations the
complex  $\M$ is connected and simply connected. 

Now, for any simple move $\mathbf{E}$ connecting parameterizations
$f,g$, we will define a functorial isomorphism $E\colon \tau [
\tsi,\{W_a\},f]\to \tau [\tsi,\{W_a\},g]$ as follows.  Recall from
\cite{TA} that simple moves are $\mathbf{Z}$
(rotation), $\mathbf{B}$ (braiding), $\mathbf{F}$ (fusion, or
gluing), $\mathbf{P}$ (which is related to isomorphsim $\phi_x$ between
standard blocks), $\mathbf{T}$ (change of marked points on a cut). 

\begin{description}
\item[$\mathbf{Z}$ move]

  If $f$ is the parameterization
  $$
  f\colon  \tsi \to S_n(g_1,g_2,\dots, g_n;h_1,h_2,\dots ,h_n)
  $$
  then $\mathbf{Z}(f)$ is the parameterization
  $$
  \mathbf{Z}(f)\colon  \tilde{\Sigma} \to
  S_n(g_n,g_1,\dots, g_{n-1};h_n,h_1,\dots ,h_{n-1})
  $$

  Now the corresponding map between $\tau [ \tsi,\{W_a\},f]$ and
  $\tau [ \tsi,\{W_a\},\mathbf{Z}(f)]$ is given by the rotation
  isomorphism of $G$-equivariant MS data which we also denoted by $Z$:
  $$
  Z\colon \< X_1,\dots,X_n  \> \to \<X_n, X_1,\dots, X_{n-1}  \>
  $$

\item[$\mathbf{B}$ move]
  Now let $f$ be the parameterization
$$
f\colon \tsi \to S_3(g_1,g_2,g_3;h_1,h_2,h_3)
$$
then $\mathbf{B}(f)$ is the parameterization
$$
\mathbf{B}(f)\colon \tilde{\Sigma} \to
S_3(g_1,g_2g_3g_2^{-1},g_2;h_1,h_3g_2^{-1},h_2)
$$
 Then the corresponding map between $\tau [ \tsi,\{W_a\},f]$
and $\tau [ \tsi,\{W_a\},\mathbf{B}(f)]$ is given by the
commutativity isomorphism of $G$ equivariant MS data which we denoted by
$\sigma$:
 $$
 \sigma \colon  \<
{}^{h_{1}^{-1}}W_1,{}^{h_{2}^{-1}}W_2,{}^{h_{3}^{-1}}W_3
\> \to \<
{}^{h_{1}^{-1}}W_1,{}^{g_{2}h_{3}^{-1}}W_3,{}^{h_{2}^{-1}}W_2 \>
 $$
\begin{remark}
Since $W_2 \in \CC_{h_2g_2h_2^{-1}}$, we must have ${}^{h_2^{-1}}W_2 \in
\CC_{g_2}$. Thus from the definition of the commutativity isomorphism,
we need to twist ${}^{h_3^{-1}}W_3$ by $g_2$.
\end{remark}

\item[$\mathbf{F}$ move] Let $f$ be the parameterization 
$$
f \colon \tsi \to S_{k+1}(\mathbf{g},\mathbf{h})
\bigsqcup_{k+1,1}
S_{l+1}(\mathbf{g'},\mathbf{h'})
$$
Here we write $\mathbf{g} = g_1,\dots, g_{k+1}$, etc. for simplicity.
Assume additionally that $h_{k+1}=h'_1$.  Then $\mathbf{F}(f)$ is the
parameterization given by
$$
\mathbf{F}(f)\colon \tsi \to
S_{k+l}(\mathbf{g''},\mathbf{h''})
$$
where 
\begin{align*}
\mathbf{g''}&=(g_1,\dots, g_k, g'_2, \dots, g'_{l+1})\\
\mathbf{h''}&=(h_1,\dots, h_k, h'_2, \dots, h'_{l+1})
\end{align*}

Then the corresponding map between $\tau [ \tsi,\{W_a\},f]$
and
$\tau [ \tsi,\{W_a\},\mathbf{F}(f)]$ is given by the gluing
isomorphism of $G$ equivariant MS data which we denoted by
$\mathcal{G}$
$$
\mathcal{G}\colon \< X_1,\dots, X_k, \RR^1 \>
\otimes \<\RR^2, Y_1,\dots,Y_l\> \to
\<X_1,\dots, X_k,Y_1,\dots,Y_l\>
$$

\item[$\mathbf{P_x}$ move]

Let $x\in\G$, and let $f$ be the
following  parameterization
$$
f \colon \tsi \to S_n(g_1,\dots,g_n;h_1,\dots,h_n)
$$
then $\mathbf{P_x}(f)$ is the parameterization
$$
\mathbf{P_x}(f)\colon \tsi \to
S_n(xg_1x^{-1},xg_2x^{-1},\dots
.xg_{n}x^{-1};h_1x^{-1},h_2x^{-1},\dots ,h_{n} x^{
-1})
$$
Then the corresponding map between $\tau [ \tsi,\{W_a\},f]$
and $\tau [ \tsi,\{W_a\},\mathbf{P_x}(f)]$ is given by the
$\phi_{x}$ isomorphism as defined in $\phi$ axiom of $G$ equivariant MS
data:
$$
\phi_{x}\colon \< {}^{h_{1}^{-1}}W_1,\dots ,{}^{h_{n}^{-1}}W_n \> \to
\< {}^{xh_{1}^{-1}}W_1,\dots,{}^{xh_{n}^{-1}}W_n \>
$$

\item[$\mathbf{T}$ move]

Let $f$ be the parameterization
$$
f\colon \tilde{\Sigma} \to S_n( \dots ,x;\dots y)
\bigsqcup_{n,1}
S_m(x^{-1},\dots  ;y,\dots )
$$
where we avoid writing all $g_i$ and $h_i$ for simplicity and only write
the labels for the boundary we want to glue. Then $\mathbf{T}$ move
replaces the label $y\in G$ by another label $z\in G$:
$$
\mathbf{T}(f)\colon \tsi \to S_n(\dots  ,x;\dots, z)
\bigsqcup_{n,1} S_m(x^{-1}, \dots ;z,\dots ..)
$$
(geometrically, it means that we are changing the choice of marked
point on the corresponding cut). Then the corresponding map between
$\tau [ \tsi,\{W_a\},f]$ and $\tau [ \tsi,\{W_a\},\mathbf{T}(f)]$ is
given by the symmetry of $\RR$ under the action of $G$, which is
part of the definition of $G$-equivariant modular functor:
$$
\< \dots \dots ,{}^{y^{-1}}\RR^1 \> \otimes \<
{}^{y^{-1}}\RR^2,\dots \> \to \< \dots \dots  ,{}^{z^{-1}}\RR^1 \>
\otimes \< {}^{z^{-1}}\RR^2,\dots\>
$$
is given by identifying both ${}^{y^{-1}}\RR^1 \tbox {}^{y^{-1}}\RR^2$
and ${}^{z^{-1}}\RR^1 \tbox {}^{z^{-1}}\RR^2$ with $\RR^1 \tbox \RR^2$
\end{description}

So far we have translated all our simple moves from the language of
parameterizations of $G$-covers to the language of $G$-equivariant MS
data. Since any two parametrizations  $f$ and $g$ can be connected  by a
sequence of simple moves (the complex $\M$ is
connected), we can construct an isomorphism of the corresponding vector
spaces, $\tau [\tsi,\{W_a\},f]$ and  $\tau [\tsi,\{W_a\},g]$. Now
we need to  show that we get the same isomorphism of the vector spaces
independent of the choice of  path (composed of simple moves). Since
the complex  $\M$ is simply connected, it is
enough to show that all the basic relations or 2-cells of
$\M$ can be translated to the corresponding axioms of
$G$-equivariant MS data.

The following is the full list of all relations (2-cells) in the
complex $\M$; precise statements of the relations can be
found in \cite{TA}.

\begin{itemize}
\item{$\mathbf{P_x}$ relation}
\item{$\mathbf{P-F}$ relation}
\item{$\mathbf{Z}$ relation}
\item{$\mathbf{B}$ relation}
\item{$\mathbf{T}$ relation}
\item{Rotation axiom}
\item{Commutativity of disjoint union}
\item{Symmetry of $\mathbf{F}$ move}
\item{Associativity of cuts}
\item{Cylinder axiom}
\item{Braiding axiom}
\item{Dehn twist axiom}
\end{itemize}

We now show how each of these relations follows from the axioms of
MS-data.

\begin{description}
\item[$\mathbf{P_x}$ relation]

  The relation $\mathbf{P_xZ = ZP_x}$
  corresponds to the compatibility of $\phi$ with rotation axiom.

    The relation $\mathbf{P_xB = BP_x}$ corresponds to the compatibility
  of $\phi$ with commutativity isomorphism.

    The relation $\mathbf{P_xF_{c,y} = F_{c,yx^{-1}}(P_x \sqcup P_x)}$
  corresponds to the compatibility of $\phi$ with gluing isomorphism.

    The relation $\mathbf{P_xP_y = P_{xy}}$ corresponds to the $\phi$
relation of the definition of MS data. 

\item[Rotation axiom]
  The rotation axiom, $\mathbf{Z}^n = \id$, of the
  complex $\M$ corresponds to the rotations axiom
  $Z^n  = \id$ of the MS data

\item[Symmetry of $\mathbf{F}$ move]
  Symmetry of the $\mathbf{F}$ move in
  $\M$ corresponds to the symmetry of gluing
  isomorphism, $\G$, in MS data

\item[Associativity of cuts]
  Associativity of cuts of
  $\M$ corresponds to the associativity of gluing
  isomorphism, $\G$, in MS data

\item[Braiding axiom]
  Braiding axioms of $\M$
  corresponds to the hexagon axiom in MS data

\item[Dehn twist axiom] 
   Dehn twist axioms of $\M$
  corresponds to the Dehn twist axiom in MS data.

\item[Cylinder axiom] 
   Cylinder axiom of $\M$
  follows from functoriality of all the isomorphisms in the definition
  of MS data.

\end{description}
We leave it to the reader to supply the details of the above
construction.

Thus, we have defined, for any   $G$-cover $\tsi$ and a parametrization
$f$,  a vector space $\tau [ \tilde{\Sigma},\{W_a\},f]$ and have shown
that for  any two prametrizations, there is a canonical isomorphism
between the vector spaces
$\tau [ \tsi,\{W_a\},f]$ and $\tau [
\tsi,\{W_a\},g]$. Now the same arguments as in $G=\{1\}$ case (see
\cite{BK}) show that this allows us to define a vector space $\tau[
\tsi,\{W_a\}]$, independent of parametrization, thus giving us the
$G$-equivaraint genus zero  modular functor.  Readers can easily show
that the $G$ equvariant genus zero modular functor defined this way
satisfies all the axioms of modular functor.

\subsection{From $G$-MF to $G$-MS data}
We assume that we are given a $\CC$ extended genus zero modular functor.

We want to create a $G$ equivarant MS data on $\CC$. To do this we
define
the following
\begin{description}
  \item[Conformal blocks]
     Given $n \geq 0$ and $W_1 \in \CC_{m_{1}}\dots,W_n \in
    \CC_{m_{n}}$ satisfying $m_1\dots m_n = 1$,
    we define
    \begin{equation}\label{e:conf_block_from_MS}
    \< W_1,\dots, W_n \> = \tau \left[
    S_{n}(m_{1},\dots,m_{n};1,\dots, 1;[W_1,\dots, W_n ]) \right]
    \end{equation}
    where $S_{n}(m_{1},\dots,m_{n};1,\dots, 1)$ is the standard block
    defined in \deref{d:stblock}.

  \item[$\phi$ axiom]
    For each $g \in G$ we define a functorial isomorphism
    $\phi_{g}\colon \< W_1,\dots, W_n \> \to 
        \< {}^gW_1,\dots, {}^g W_n \>$ as the following
    composition:
    \begin{align*}
    &\tau \left[S_{n}(m_{1},\dots,m_{n};1,\dots,1;[W_1,\dots, W_n])
        \right] \\
    &\qquad \xxto{(\phi_g)_*}
    \tau\left[
      S_{n}(gm_{1}g^{-1},\dots,gm_{n}g^{-1};g^{-1},\dots,g^{-1};
            [W_1,\dots,W_n]) \right]\\
    &\qquad \xxto{T_{g,\dots,g}}
    \tau \left[
          S_{n}(gm_{1}g^{-1},\dots,gm_{n}g^{-1};1,\dots,1;
                [{}^gW_1,\dots, {}^g W_n]) 
       \right]
  \end{align*}
  Here $\phi_g$ is the  morphism between standard blocks described in 
  \leref{isolemma}, and $T_{g}$ is as in the definition of the modular
  functor (see eq. \eqref{e:Tx}). 
%

  \item[Rotation axiom]
    The rotation isomorphism
    $$
    Z\colon \< W_1,\dots,W_n \> \to 
            \< W_n, W_1,\dots,W_{n-1} \>
    $$
    is given by $Z=(\tilde z)_*$, where $\tilde z$ is the rotation 
    homeomorphism of the standard block (see eq. \eqref{e:z}). 

  \item[Symmetric object]
    The  symmetric object $\RR$ directly comes from the definition of
    $G$ equivariant modular functor. 

  \item[Gluing isomorphism]
    Let $A_{i} \in \CC_{p_{i}}$ and $B_{j} \in
    \CC_{q_{j}}$, where $i=1\dots k$ and $j=1\dots l$ and they satisfy
    $p_1\dots p_kq_1\dots q_l=1$. Then the gluing isomorphism
     $$
     \mathcal{G}\colon \< A_1,\dots,A_k,\RR^1 \> \otimes 
                       \< \RR^2,B_1,\dots,B_l \> 
            \rightarrow \< A_1,\dots,A_k,B_1,\dots,B_l \>
     $$
     is given by the $(\tilde{\al}_{k,l})_*$, where 
     \begin{align*}
     \tilde{\al}_{k,l}&\colon 
           S_{k+1}(p_{1},\dots,p_{k},x;1,\dots,1)
           \bigsqcup_{k, 1}
           S_{l+1}(x^{-1},q_{1},\dots,q_{l};1,\dots,1)\\
      &\rightarrow
     S_{k+l}(p_1,\dots,p_k,q_1,\dots,q_l;1,\dots,1)
     \end{align*}
     is the homeomomorphism defined in  paper \cite{TA}.
     Note that we must have  $x=q_1\dots q_l=(p_1\dots
      p_k)^{-1}$. 

  \item[Commutativity isomorphism]
   The commutativity isomorphisms 
     \begin{align*}
     &\sigma\colon  \< X, A,B \> \to \< X,  {}^{p}B,A\>\\
      &X \in \CC_r,\quad A \in \CC_p,\quad B\in \CC_q, \quad rpq= 1
     \end{align*}
      is defined as follows. Recall the braiding homeomorphism (see
      \cite{TA}):
     $$
     \tilde{b}\colon  S_3(r, p,q;1,1,1) \to
S_3(r, pqp^{-1},p;1, p^{-1}, 1)
     $$
     Then we define $\si$ as the following composition: 
    \begin{align*}
    &\tau [ S_3(r, p,q;1,1,1;X, A,B)] \xxto{\tilde b_*}  
     \tau [S_3(r, pqp^{-1},p; 1, p^{-1},1;X, B,A)]\\
     &\qquad  \xxto{T_{1, p,1}} 
               \tau [S_3(r, pqp^{-1},p;1,1,1;X, {}^pB,A)]
    \end{align*}
 
\end{description}

This completes the construction of MS data from a MF.

Now we need to check that so defined isomorphisms satisfy all the axioms
of $G$-equivariant Moore-Seiberg data.
%
%
%
%
%
\begin{itemize}
  \item 
      Non-degeneracy of MS data follows from the non-degeneracy of
      modular functor.
  \item 
      Normalization axiom of MS data follows from the normalization
      axiom of modular functor.
  \item 
    Associativity of $\mathcal{G}$ of MS data follows from the
    associativity of cuts of the complex $\M$
  \item 
    Rotation axiom of MS data follows from the rotation axiom of the
    complex $M(\tsi,\Sigma)$.
  \item 
    Symmetry of $\mathcal{G}$ in MS data follows from the symmetry
    of $\mathbf{F}$ move in the complex $\M$.

  \item 
    $\phi$-relation and the compatibility of $\phi$ follows from the
    $\mathbf{P_x}$ relation of the complex $\M$
    and the associativity of the group multiplication in $G$.

  \item 
    Hexagon axiom of MS data follows from the braiding axiom of the
    complex $\M$. 

  \item 
    Dehn twist axiom of MS data follows from the Dehn twist axiom of
    the complex $\M$. 
\end{itemize}
Most of the above correspondences are  easy to check since they
directly follow from the definition. For illustration, we will
demonstrate   the Dehn twist axiom  for  MS data.

If we look at the Dehn twist axiom of MS data and rewrite everything
using our definition of conformal block, replacing  $\<A,B\>$ by 
$\tau[S_2(p,q;1,1;[A,B])]$ etc, we get the following diagram (as
before, $A\in\CC_p, B\in\CC_q, pq=1$. 
      \begin{figure}[ht]
      $$ 
       \xymatrix@R=5pt{
       \tau[S_2(p,q;1,1;[A,B])] \ar[r]^\si \ar[dd]_Z & 
                       \tau[S_2(q,p;1,1;[{}^pB,A])]\ar[dr]_Z\\
        & & \tau[S_2(p,q;1,1;[A,{}^pB])]\\
       \tau[S_2(q,p;1,1;[B,A])]\ar[r]^\si 
           &\tau[S_2(p,q;1,1;[{}^qA,B])]\ar[ur]_\ph  
       }$$

     \caption{Rewriting Dehn twist axiom}\label{f:dehn_twist_explained}
      \end{figure}

Now it is  easy to derive the Dehn twist
axiom of MS data from the axioms of modular functor. Recall that in 
the complex 
$\M$, we had the following relation,
which we also called ``Dehn twist axiom'':  (we assume
$\alpha$ is the boundary associated with $A$ and the label $(p,1)$ and
$\beta$ is the boundary circle associated with $B$ and the label
$(q,1)$; as before,  $pq = 1$)
 $$
  \xymatrix@R=5pt{
       S_2(p,q;1,1)\ar[dd]_Z \ar[r]^{B_{\alpha, \beta}} & 
                       S_2(q,p;q,1)\ar[dr]^Z\\
    &&S_2(p,q;1,q)\\
     S_2(q,p;1,1)\ar[r]^{B_{\be,\al}} &S_2(p,q;p,1)\ar[ur]^{\ph_p}
}
$$
The above relation is a certian relation between homeomorphsism of
$G$-covers; bu functoriality axiom of modular functor, this relation
must also hold between the corresponding modualr functor spaces, which
exactly gives us realtion of \firef{f:dehn_twist_explained}, this
proving the Dehn twist axiom for MS data.

The proof of all othe axioms of MS data is done in a similar manner, by
rewriting the axioms in terms of modular functro spaces for standard
blocks and using relations in the complex $\M$.

This finishes the proof of our main theorem.



\end{document}